\pdfoutput=1
\documentclass{amsart}
\usepackage{colortbl}
\usepackage{amsmath,amsthm,amssymb}
\usepackage{graphicx}
\usepackage{harvard}

\numberwithin{equation}{section}

\def\ZZZ{{\mathbb Z}}
\def\LLL{{\mathcal L}}

\def\RRR{{\mathbb R}}

\def\RRR{{\Bbb R}}

\def\1{\mathbf 1}

\newtheorem{definition}{Definition}
\newtheorem{notation}{Notation}
\newtheorem{theorem}{Theorem}
\newtheorem{lemma}{Lemma}
\newtheorem{corollary}{Corollary}

\begin{document}

\title{Stochastic Integrals and Abelian Processes}

\author{Claudio Albanese}

\email{claudio@level3finance.com}

\date{First version November 1st, 2007, last revision \today}

\thanks{}

\maketitle

\begin{abstract}

We study triangulation schemes for the joint kernel of a diffusion
process with uniformly continuous coefficients and an adapted,
non-resonant Abelian process. The prototypical example of Abelian
process to which our methods apply is given by stochastic integrals
with uniformly continuous coefficients. The range of applicability
includes also a broader class of processes of practical
relevance, such as the sup process and certain discrete time
summations we discuss.

We discretize the space coordinate in uniform steps and assume that time is either
continuous or finely discretized as in a fully explicit Euler
method and the Courant condition is satisfied.
We show that the Fourier transform of the joint kernel of a
diffusion and a stochastic integral converges in a uniform graph
norm associated to the Markov generator. Convergence also implies
smoothness properties for the Fourier transform of the joint kernel.
Stochastic integrals are straightforward to define for finite
triangulations and the convergence result gives a new and entirely
constructive way of defining stochastic integrals in the continuum.
The method relies on a reinterpretation and extension of the classic
theorems by Feynman-Kac, Girsanov, Ito and Cameron-Martin, which are
also reobtained.

We make use of a path-wise analysis without relying on a
probabilistic interpretation. The Fourier representation
is needed to regularize the hypo-elliptic character of the joint
process of a diffusion and an adapted stochastic integral.
The argument extends as long as the Fourier analysis framework can
be generalized. This condition leads to the notion of non-resonant
Abelian process.

\end{abstract}

\tableofcontents

\section{Introduction}

Consider a diffusion defined by the Markov generator
\begin{equation}
{\mathcal L}^{0}_x = {1 \over 2} {\sigma(x)^2 } {\partial^2\over
\partial x^2} + \mu(x){\partial \over \partial x}.
\end{equation}
on the bounded interval $X = [-L_{x},L_{x}] \subset \RRR$ where
$0<L_{x}<\infty$. For simplicity, we assume periodic boundary
conditions and identify the two boundary points $\pm L_{x}$ with
each other. The coefficients $\sigma(x)^2$ and $\mu(x)$ are assumed
to be at least uniformly continuous and our results will depend on
the degree of smoothness. We also assume that $\sigma(x)>\Sigma_0$
for some constant $\Sigma_0>0$.

In this paper, we consider path dependent processes such as
stochastic integrals of the form
\begin{equation}
y_t = \int_0^t a(x_s) d x_s + b(x_s) ds \label{stochint}
\end{equation}
where $a(x)$ and $b(x)$ are functions on $A = [-L_{x}, L_{x}]$. The
case $a(x) = 0$ was considered in \cite{Girsanov},
\cite{CameronMartin}, \cite{Feynman1948} and \cite{Kac1950}. The
case with $a(x) \neq 0$ was first tackled in \cite{Ito1951}. The
proofs were streamlined by using martingale theory in
\cite{Doob1953}. Other references for stochastic integrals are
\cite{KunitaWatanabe1967}, \cite{McK} and \cite{Meyer1974}. Most of
these results assume smooth coefficients. The more general case of
uniformly continuous coefficients was considered in \cite{SV1969}.

We are interested in convergence properties of triangulation schemes
for the joint distribution between the underlying process $x_t$ and
the stochastic integral $y_t$. If the diffusion process is
approximated by a sequence of Markov chains then existence and
uniqueness of stochastic integrals is not problematic. The problem
is to determine conditions under which joint distributions between
stochastic integrals and the underlying process converge in some
meaningful norm in the limit as the triangulation becomes finer and
finer, approaching the continuum limit. We establish convergence
assuming uniform continuity in the coefficients of the diffusion and
the stochastic integral and find convergence rates. Convergence
takes place in the graph-norm of the Fourier transformed Markov
generator. In the continuum limit, the Fourier transform of the
joint distribution is entire analytic in the conjugate variable of
the stochastic integral and is in the domain of all the powers of
the Fourier transformed Markov generator. Consideration of the
Fourier transform is essential in the proof as the joint
distribution itself can be singular even in trivial special cases
such as the one where $a(x) = b(x) = 0$. The Fourier analysis also
leads to a new derivation of the celebrated formulas for
characteristic functions of stochastic integrals in \cite{Girsanov},
\cite{CameronMartin}, \cite{Feynman1948} and \cite{Ito1951}.
Furthermore, analyticity implies convergence of moment formulas.

These results are classic and well established in the continuum
limit. One of the new facets of our derivation is that it is
entirely constructive and thus instructive from a computational
viewpoint. We make no use of compactness arguments and
non-constructive measure theory methods. The derivation is based
on pathwise analysis and a renormalization group transformation.
The construction applies to discretized operators and does not rely
on the existence of a probabilistic interpretation. The proof of convergence
is actually carried out for analytic extensions of the Fourier transform
of the joint distribution, so we consider complex valued weights for
the paths. The argument would also extend
to different type of diffusions with complex coefficients as they
occur for instance in quantum mechanics.
Furthermore results are extended to discrete time approximations
with fully explicit Euler schemes, in which case we also derive
convergence bounds in the graph-norm for Fourier transformed joint
kernels.

In \cite{AlbaneseKernelEsimtatesA} we discuss the case of the kernel
of diffusion processes, a more elementary situation. This paper is
strictly more general and it is self-contained with no dependencies
on our previous work as all derivations need to be adapted in
detail. The key technique of pathwise analysis and renormalization
group resummations over decorated paths are however quite similar.

Our operator algebraic approach to the problem also points out that
stochastic integrals are not the most general class of path
dependent processes to which these methods apply. This is based on
the recognition that Fourier transforms play a key role as they
allow one to block-diagonalize the joint kernel. But
block-diagonalization can be achieved under more general conditions.
We consider a generic situation of a diffusion and an adapted
process. We associate an operator algebra to a finite triangulation of
the joint process. If this operator algebra is Abelian, i.e.
commutative, and if it satisfies a certain non-resonance condition,
then the joint Markov generator can also be block-diagonalized and
methods can be extended.

Non-resonant Abelian processes are a broad class of hypo-elliptic
multidimensional diffusion processes which are computable by robust
methodologies based on block-diagonalizations and fully explicit
Euler schemes. Several instances of these processes found
engineering applications already, see \cite{AVidler}, \cite{ATrovato2},
\cite{AlbaneseVideogames2007}, \cite{ALoMijatovic}, \cite{AOsseiran}.
As an example, we discuss here the case of the sup processes of a one-dimensional
diffusion, i.e. of
\begin{equation}
y_t = \sup_{\in [0,t]} x_s.
\label{supprocess}
\end{equation}
We also discuss discrete time summations of the form
\begin{equation}
y_{T_n} = \sum_{i=1}^n \phi(x_{T_{i-1}}, x_{T_i})
\label{supprocess}
\end{equation}
where $T_i = i \Delta T$.

The kernel convergence in graph norm we establish in this paper has
noteworthy implications from the computer science viewpoint. As we
discuss in \cite{AlbaneseFundamental2007}, they imply that explicit
schemes give a robust valuation methodology for joint kernels of
processes specified semi-parametrically. The Courant-Friedrichs-Lewi
stability condition \cite{CFL} for explicit methods is not much of an impediment
in the presence of sufficient system memory as linear fast exponentiation
can be used to accelerate the scheme.
The ability of evaluating kernels robustly is
the main reason why analytically solvable models are interesting and
broadly used. But not all models solvable in closed form are
practically computable in the sense that they allow for robust
kernel valuations within the limits of double precision floating
point arithmetics. Higher
order hypergeometric functions are often difficult to handle and
may require multiple precision. Linear
fast exponentiation is a technique based on accelerated fully
explicit schemes that empirically is observed to work well even in
single precision and for very general model specifications. Empirically, we
find that it performs better in single precision on some
solvable models than the
use of closed form solutions in double precision. This paper was
prompted by the observation of this phenomenon and the desire to
understand and explain the underlying smoothing mechanisms which is
of technological significance.

The paper is organized as follows: In Section \ref{SecStochint} we
introduce notations and state our results regarding stochastic
integrals. Proofs in the case of stochastic integrals are in Section
\ref{SecRenormalization}, where we consider the case of continuous
time and in Section \ref{SecEuler} where we discuss convergence for
fully explicit Euler schemes. The Ito representation is in Section 5.
Extensions to the sup and other Abelian processes are given in Section 6
and conclusions end the paper.

\section{Stochastic Integrals \label{SecStochint}}

Consider a family of increasingly fine triangulations schemes
whereby the space dimension is discretized in multiples of an
elementary space step $h_{xm} = L_{x} 2^{-m}, m \in {\Bbb N}$ and we
are interested in the limit as $m\to\infty$. Let $X_m =
{h_{xm}}\ZZZ\cap X$ and consider the sequence of operators
\begin{equation}
{\mathcal L}^{m}_x = {\sigma(x)^2 \over 2} \Delta^{m}_x +
\mu(x)\nabla^{m}_x. \label{eq_lm}
\end{equation}
defined on the $2^{m+1}$-dimensional space of all periodic functions
$f_{m}: X_m \to \RRR$, where
\begin{equation}
\nabla^{m}_x f(x) = {f(x+{h_{xm}})-f(x-{h_{xm}})\over 2 {h_{xm}}}.
\end{equation}
and
\begin{equation}
\Delta^{m}_x f(x) = {f(x+{h_{xm}})+f(x-{h_{xm}})-2 f(x) \over
{h_{xm}}^2}
\end{equation}
These definitions also apply to the boundary points by periodicity.
We assume that $m\ge m_0$ where $m_0$ is the least integer such that
\begin{equation}
{\sigma^2(x)\over 2 h_{xm}^2} > {\lvert \mu(x) \lvert \over 2h_{xm}}
\label{def_diffusion_bounds}
\end{equation}
for all $m\ge m_0$ and all $x\in X_m$.

Consider the kernel $u_{m}(x, x'; t)$ of equation (\ref{eqdiff}),
i.e. the solution of the (forward) equation
\begin{equation}
{\partial \over \partial t} u_{m}(x, x'; t) = {\mathcal L}^{m
*}_{x'} u_{m}(x, x'; t) \label{eqdiff}
\end{equation}
where the operator ${\mathcal L}^{m *}_{x'}$ acts on the $x'$
coordinate and the following initial time condition is satisfied:
\begin{equation}
u_{m}(x, x'; 0) = h_{xm}^{-1} \delta_{X_m}(x-x'). \label{eqdiff}
\end{equation}
where
\begin{equation}
\delta_{X_m}(x-x') =
\begin{cases}
1 \;\;\;\;\;\;\;\;\;{\rm if} \;\; x = x' \;\;\;{\rm mod} \;2 L_{x}\\
0 \;\;\;\;\;\;\;\;\;{\rm otherwise}.
\end{cases}
\label{eqdelta}
\end{equation}
If $f(x, x'):X_m\times X_m \to \RRR$, consider the uniform norm
\begin{equation}
\lvert\lvert f_m \lvert\lvert_{m,\infty} = \sup_{x, x' \in X_m}
\lvert f(x, x') \lvert
\end{equation}
and the graph norm
\begin{equation}
\lvert\lvert f \lvert\lvert_{m, \mathcal L} = \lvert\lvert f
\lvert\lvert_\infty + \lvert\lvert {\mathcal L}^m_x f
\lvert\lvert_\infty + \lvert\lvert {\mathcal L}^{m*}_{x'} f
\lvert\lvert_\infty.
\end{equation}
In \cite{AlbaneseKernelEsimtatesA} we show that if the coefficients
$\sigma(x)^2$ and $\mu(x)$  are uniformly continuous, then the
sequence of kernels $u_m$ is Cauchy with respect to the graph norm.
If in addition coefficients are H\"older continuous and $\sigma^2
\in {\mathcal C}^{k_\sigma, \alpha_\sigma}$, $\mu \in {\mathcal
C}^{k_\mu, \alpha_\mu}$ and
\begin{equation}
\gamma \equiv \min\{ 2, k_\sigma+\alpha_\sigma, k_\mu+\alpha_\mu \}
>0 \label{eq_gamma}
\end{equation}
then
\begin{equation} \lvert\lvert u_m(t) - u_{m'}(t)
\lvert\lvert_{m, \mathcal L} \leq c h_{xm}^\gamma
\end{equation}
for all $m'>m\ge m_0$. A similar bound also holds for the kernels
obtained with a fully explicit Euler scheme
\begin{equation}
u^{\delta t}_m(x, x'; t) = h_{xm}^{-1} \left( 1 + \delta t_m
{\mathcal L}^m \right)^{\left[t\over\delta t_m\right]} (x, x'; t).
\end{equation}
where $\delta t_m $ is so small that
\begin{equation}
\min_{x\in X_m} 1 + \delta t_m {\mathcal L}^m(x, x) > 0.
\end{equation}
In fact, we have that
\begin{equation} \lvert\lvert u_m(t) - u^{\delta
t}_{m'}(t) \lvert\lvert_{m, \mathcal L} \leq c h_{xm}^2
\end{equation}
for some constant $c>0$.

Let $h_{yn}$ be a monotonously descreasing sequence such as $h_{yn}
\to 0$ as $n\to\infty$. Also let $L_{yn}$ be a monotonously
increasing sequence such that $L_{yn}\to\infty$ as $n\to\infty$. Let
\begin{equation}
Y_n = \big(h_{yn} \ZZZ\big) \cap [-L_{yn}, L_{yn}].
\end{equation}
We assume again periodic boundary conditions in the $y$ direction
and identify the two extreme points of $Y_n$. Let's consider
two-dimensional processes described by a sequence of Markov
generators of the form
\begin{equation}
{\mathcal L}^{m, n}  (x, y; x', y') =  {\mathcal L}^m (x, x') +
{\mathcal Q}^{m, n} (x, y; x', y') \label{eq_ljoint}
\end{equation}
where $x, x' \in X_n$ and $y, y' \in Y_n$. We assume that
\begin{equation}
\sum_{y'\in Y_n} {\mathcal Q}^{m, n} (x, y; x', y') = 0
\end{equation}
for all values of $x, y, x'$, so that the marginals with respect to
the first process are the same as under the dynamics given by
${\mathcal L}^{m}$, i.e.
\begin{equation}
\sum_{y'\in Y_n} \exp\big( t {\mathcal L}^{m, n}\big)  (x, y; x',
y') = \exp\big( t {\mathcal L}^{m}\big)  (x; x')
\end{equation}
for all triples $x, y, x'$. We are interested in the joint kernel
\begin{equation}
u_{mn}(x, y; x', y'; t) \equiv {1\over h_{xm} h_{yn}} \exp\big(t
{\mathcal L}^{m, n}\big)(x, y; x', y')
\end{equation}
and its convergence properties in the limit as $m, n \to \infty$.

The case of general multi-dimensional diffusion processes is
discussed in \cite{AlbaneseJones2007}. Results in that paper however
do not apply to the class of processes in which we are interested
here. In this article, we consider the case of path dependent
processes which are not driven by their own diffusive dynamics and
whose generator is not strongly elliptic, but rather hypoelliptic.
Because of this reason, the joint kernel may be singular and a
separate treatment is required.

\begin{definition}
Consider a sequence of Markov generators of the form in
(\ref{eq_ljoint}) describing a bivariate process. One says that the
second process is Abelian with respect to the first if the matrices
$Q^{m, n}_{x, x'}$ of elements
\begin{equation}
Q^{m, n}_{x, x'} (y, y') = {\mathcal Q}^{m, n} (x, y; x', y')
\end{equation}
are mutually commuting, i.e. if
\begin{equation}
[Q^{m, n}_{x, x'} , Q^{m, n}_{x'', x'''}] = Q^{m, n}_{x, x'} Q^{m,
n}_{x'', x'''} - Q^{m, n}_{x'', x'''}  Q^{m, n}_{x, x'} = 0
\end{equation}
for all $x, x', x'', x''' \in X_m$.
\end{definition}

\begin{definition}
An Abelian process is called non-resonant if there is a sequence of
non-singular transformations ${\mathcal V}_{mn}(p, y)$ which
diagonalizes simultaneously the kernels $Q^{m, n}_{x, x'}$ for all
pairs $x, x' \in X_m$, i.e. if
\begin{equation}
\sum_{y, y' \in Y_{n}} {\mathcal V}_{mn}(p, y) {\mathcal Q}^{m, n}
(x, y; x', y') {\mathcal V}_{mn}^{-1}(p', y') = \Lambda_{mn}(x, x')
\delta_{p, p'}.
\end{equation}
for all $x, x' \in X_m$. The index $p$ ranges on a set denoted with
$\hat Y_{n}$ which has the same cardinality as $Y_{n}$ and is called
the inverse lattice of $Y_{n}$.
\end{definition}

In Sections \ref{SecRenormalization}, \ref{SecEuler}
and \ref{SecItoRepresentation}, we consider in detail the prototypical
example of Abelian process, namely stochastic integrals. In this case,
a triangulation for the joint generator is given by
\begin{align}
{\mathcal L}^{m, n}  (x, y; x', y') =& \left( {\sigma(x)^2\over 2
h_{xm}^2} + {\mu(x)\over 2 h_{xm}} \right) \delta_{X_m}(x' - x -
h_{xm}) \delta_{Y_n}\bigg(y' - y
-\bigg[{a(x) h_{xm} \over h_{yn}}\bigg] h_{yn} \bigg) \notag  \\
&+ \left( {\sigma(x)^2\over 2 h_{xm}^2} - {\mu(x)\over 2 h_{xm}}
\right) \delta_{X_m}(x' - x + h_{xm}) \delta_{Y_n}\bigg(y' - y
+\bigg[{a(x) h_{xm}
\over h_{yn}}\bigg] h_{yn} \bigg) \notag \\
&- {\sigma(x)^2\over h_{xm}^2} \delta_{X_m}(x' - x) \delta_{Y_n}(y'
- y)
+   {b(x)} \delta_{x x'} \nabla^n_y (y, y') \notag \\
\label{eq_stochintl}
\end{align}
where
\begin{equation}
\nabla^{n}_y g(y) = {g(y+{h_{yn}})-g(y-{h_{yn}})\over 2 {h_{yn}}}.
\end{equation}
and
\begin{equation}
\delta_{Y_n}(y-y') =
\begin{cases}
1 \;\;\;\;\;\;\;\;\;{\rm if} \;\; y = y' \;\;\;{\rm mod} \;2 L_{y}\\
0 \;\;\;\;\;\;\;\;\;{\rm otherwise}.
\end{cases}
\label{eqdeltay}
\end{equation}
Diagonalisation is simply achieved by Fourier transforms of kernel
\begin{align}
{\mathcal F}_{n}(x, p; x', y) = e^{-i y p} \delta_{x x'}
\end{align}
where
\begin{equation}
p\in \hat Y_n \equiv \big(\hat h_{yn}\ZZZ\big) \cap \left[ - {\pi
\over h_y} , {\pi \over h_y} \right), \;\;\;\; {\rm and}
\;\;\;\;\hat h_y = {\pi\over L_{yn}}.
\end{equation}
Partial Fourier transforms in the $p$ variable block-diagonalize the
joint generator, reducing it to the form
\begin{align}
\big({\mathcal F}_{n} {\mathcal L}^{m, n}{\mathcal F}_{n}^{-1} \big)
(x, p; x', p') = \hat {\mathcal L}^{m, n}(x, x'; p) \delta_{p p'}
\end{align}
where  $\hat {\mathcal L}^{m, n}(x, x'; p)$ is a one-parameter
family of matrices indexed by $x, x' \in X_m$. An explicit
expression for the joint generator is given in equation
(\ref{eq_fourierl}) below. The kernel itself can also be expressed
by means of a Fourier transform as follows:
\begin{equation}
u_{mn}(x, p; x', p'; t) \equiv \big({\mathcal F}_{n}^{-1} \hat
u_{mn}(t) {\mathcal F}_{n}\big)(x, p; x', p') = {1\over h_{xm} h_{yn}}
\delta_{\hat Y_n}(p-
p') \exp\left( t \hat {\mathcal L}^{m, n}(p) \right)(x, x') ,
\;\;\;\;\;
\end{equation}
We are interested in the convergence properties of the sequence
\begin{equation}
\hat u_{m}(x, x'; p, t) =  \lim_{n\to\infty} {1\over h_{yn}} \hat
u_{mn}(x, p; x', p'; t)
\end{equation}
in the limit as $m\to\infty$. To establish convergence, it is
necessary to work in the representation of the Fourier transformed
kernel. In fact, while the Fourier transformed kernel converges in a
very strong graph-norm and the limit is entire analytic in $p$ and
has the same smoothness properties in $x, x'$ as the diffusion
kernel, the joint kernel expressed with respect to space coordinates
can be quite singular even in simple cases. If for instance the
functions $a(x)$ and $b(x)$ in (\ref{stochint}) are zero, then the
joint distribution $U(x, y; x', y'; t)$ for the pair $(x_t, y_t)$ is
concentrated on the line $y_t=0$. If a partial Fourier transform
with respect to the $y$ variable is taken, the result is constant as
a function of the dual variable $p$. This is a simple example of
entire analytic dependency in $p$ which however translates into a
rough delta type singularity in the coordinate representation.

Let $g(x, x'; z)$ be a complex valued function defined for $x, x'\in
X_m\times X_m$ and for $z\in{\Bbb C}$ such that $\lvert z \lvert <
K$ where $K>0$. Consider the uniform norm
\begin{equation}
\lvert\lvert g \lvert\lvert_{m, K, \infty} = \sup_{\begin{matrix} x,
x' \in X_m \\ \lvert z \lvert < K \end{matrix}} \lvert g(x, x' ; z)
\lvert.
\end{equation}
The graph norm of order $(m,K,
{\mathcal L})$ is defined as follows:
\begin{equation}
\lvert\lvert g \lvert\lvert_{m , K, {\mathcal L}} =
\lvert\lvert g_{m} \lvert\lvert_{m, K, \infty} + \lvert\lvert
{\mathcal L}^{m}_{x} g \lvert\lvert_{m, K, \infty} + \lvert\lvert
{\mathcal L}^{m*}_{x'} g \lvert\lvert_{m, K, \infty}.
\end{equation}

\begin{theorem}{\label{kernel_convergence}} Let's assume the function $b(x)$ is integrable and
let
\begin{equation}
B(x) = \int_0^x b(x') dx'
\end{equation}
be a primitive. If the coefficients $\sigma(x)^2$, $\mu(x)$ and the
functions $a(x)$ and $B(x)$ are uniformly continuous, then the
sequence of Fourier transformed kernels $\hat u_m(t)$ can be
extended by analyticity in $p$ to an entire operator valued function
and is Cauchy with respect to the above graph-norm for all $K>0$. If
in addition these coefficients are H\"older continuous so that
$\sigma^2 \in {\mathcal C}^{k_\sigma, \alpha_\sigma}$, $\mu \in
{\mathcal C}^{k_\mu, \alpha_\mu}$, $a \in {\mathcal C}^{k_a,
\alpha_a}$, $B \in {\mathcal C}^{k_B, \alpha_B}$ and
\begin{equation}
\gamma \equiv \min\{ 2, k_\sigma+\alpha_\sigma, k_\mu+\alpha_\mu\ ,
k_a+\alpha_a, k_B +\alpha_B\}
>0 \label{eq_gamma}
\end{equation}
then, for all $K>0$ there is a constant $c(K)$ such that
\begin{equation} \lvert\lvert \hat u_{m}(t) - \hat u_{m'}(t)
\lvert\lvert_{m, K, {\mathcal L}} \leq c(K) h_{xm}^{\gamma}.
\end{equation}
for all $m'>m\ge m_0$.
\end{theorem}

\begin{corollary}
The limit kernel
\begin{equation}
u(x, x'; z) \equiv \lim_{m\to\infty} {1\over h_{xm}} u_m(x, x'; z)
\end{equation}
is an entire analytic function of $z$ which is in the domain of the
operators ${\mathcal L}_x$ and ${\mathcal L}_y^*$ and all of their powers.
The Fourier transformed kernel
\begin{equation}
u(x, x'; y) \equiv \int_{-\infty}^\infty e^{-i p y} u(x, x'; p) {dp \over 2\pi}
\end{equation}
exists in the distribution sense for each fixed pair $x, x' \in [-L_x, L_x]$.
\end{corollary}

\begin{theorem}{\label{kernel_euler}}
Let us assume that the coefficients $\sigma(x)^2$, $\mu(x)$ and the
functions $a(x)$ and $B(x)$ are uniformly continuous and
fix a $K>0$. Consider the kernels obtained with a fully explicit Euler scheme, i.e.
\begin{equation}
\hat u^{\delta t}_{m}(x, x'; z, t) = {1\over h_{x}} \left( 1 +
{\delta t}_{m} \hat {\mathcal L}^{m}(z) \right)^{\left[t\over\delta
t_{m}\right]}(x, x'),
\end{equation}
where $z\in{\Bbb C}: \lvert z \lvert < K $. Assume that
$\delta t_{m} $ satisfies the Courant condition
\begin{equation}
\min_{x, p\in X_m\times \hat Y_n} \Re\bigg( 1 + \delta t_{m} \hat {\mathcal
L}^{m}(x, x; z)\bigg)
> 0
\end{equation}
for all $\lvert z \lvert<K$. Then, there is a constant $c(K)$ such that
\begin{equation} \lvert\lvert \hat u_{m}(t) - \hat u^{\delta t}_{m}(t)
\lvert\lvert_{m, K, {\mathcal L}} \leq c(K) h_{xm}^2
\end{equation}
for all $m'>m\ge m_0$.
\end{theorem}

This two Theorems are proved in the next two sections.
An equivalent statement in a different
representation which leads to Ito's Lemma, the Cameron-Marin-Girsanov
and the Feynman-Kac formulas is given in Section \ref{SecItoRepresentation}.

\section{The Renormalization Group Argument \label{SecRenormalization}}

The proof of convergence in graph-norm is based on a path-integral
representation of the probability kernel. More precisely, let us
defines a {\it symbolic path} $\gamma = \{\gamma_0, \gamma_1,
\gamma_2, .... \}$ as an infinite sequence of sites in $X_m$ such
that $\gamma_j = \gamma_{j-1} \pm 1$ for all $j=1, ...$. Let
$\Gamma_m$ be the set of all symbolic paths in $X_m$. Then the
propagator admits the following representation:
\begin{align}
u_m(x, x'; t) =& {1\over h_{xm}} \sum_{q=1}^\infty \sum_{\gamma\in
\Gamma_m : \gamma_0 = x, \gamma_{q} = x'}
 \int_0^{t} ds_1 \int_{s_1}^{t} ds_2 ...
\int_{s_{q-1}}^{t} d s_{q} \rho(s_1, s_2, ..., s_{q-1} ; \gamma)
\label{eq_upathint}
\end{align}
where
\begin{align}
\rho(s_1, s_2, ..., s_{q-1} ; \gamma)=& e^{-s_{1}
\LLL^m(\gamma_0, \gamma_0) } \prod_{j=1}^{q-1} \LLL^m(\gamma_{j-1},
\gamma_j) e^{ - (s_{j+1}- s_{j}) \LLL^m(\gamma_j, \gamma_j)} ds_1
... ds_{q}. \label{eq_upathint}
\end{align}
and $s_{q} = t$.

Let $\gamma\in\Gamma_m$ and consider a path $X_m(\cdot; t_1, ...
t_q; \gamma): \RRR_+ \to X_m$, left continuous and with right
limits, taking the values $\{\gamma_0, \gamma_1, \gamma_2, .... \}$
consecutively with jumps occurring at times $t_j, j=1,.. q-1$, with
$0\leq t_1 \leq t_2 \leq .. \leq t_{q-1}\leq t$. Let $a(x)$ and
$b(x)$ be two uniformly continuous functions in $X = [-L_x, L_x]$
and consider the integral
\begin{align}
I(t_1, ... t_{q-1}; \gamma) & = \int_0^t \bigg[ a(X_m(s; t_1, ...
t_q; \gamma)) { d X_m(s; t_1, ... t_q; \gamma) \over ds}
 + b(X_m(s; t_1, ... t_q;
\gamma)) \bigg] ds \notag \\
&\hskip5cm=\sum_{j=0}^{q-1} a(\gamma_{j}) (\gamma_{j+1} -
\gamma_{j}) +
 \sum_{j=0}^q b(\gamma_{j}) ({t_{j+1}} - {t_{j}}). \notag \\
\end{align}
We are interested in evaluating the joint density
\begin{align}
u_m(x, y; x', y'; t) =& {1\over h_{xm}} \sum_{q=1}^\infty \sum_{\gamma\in
\Gamma_m : \gamma_0 = x, \gamma_{q} = x'}
 \int_0^{t} ds_1 \int_{s_1}^{t} ds_2 ...
\int_{s_{q-1}}^{t} d s_{q} \notag \\
&\hskip5cm \rho(s_1, s_2, ..., s_{q-1} ; \gamma)
\delta\big(I(s_1, ... s_{q-1}; \gamma) - (y'-y)\big) \notag \\
\label{eq_upathint}
\end{align}
We do so by expanding its Fourier transform
\begin{align}
\hat u_m(x, x'; p; t) =& {1\over h_{xm}} \sum_{q=1}^\infty
\sum_{\gamma\in \Gamma_m : \gamma_0 = x, \gamma_{q} = x'}
 \int_0^{t} ds_1 \int_{s_1}^{t} ds_2 ...
\int_{s_{q-1}}^{t} d s_{q} \rho(s_1, s_2, ..., s_{q-1} ; \gamma)
e^{i p I(s_1, ... s_{q-1}; \gamma)} \label{eq_upathint}
\end{align}

The partial Fourier transform of the joint generator in
(\ref{eq_stochintl}) in the limit as $n\to\infty$ is given by
\begin{align}
\hat\LLL^{m}(x, x'; p) & =  \lim_{n\to\infty}
\hat\LLL^{mn}(x, x'; p) \notag \\
& = \lim_{n\to\infty}
 \sum_{y\in Y_n} \tilde\LLL^{mn}(x, 0; x', y; t)
e^{- i h_{yn} p y } \notag = \left( {\sigma(x)^2\over 2 h_{xm}^2} +
{\mu(x)\over 2 h_{xm}} \right) e^{- i {h_{xm} a(x)} p } \delta_{x',
x+h_{xm}} \notag \\
& \hskip3cm + \left( {\sigma(x)^2\over 2 h_{xm}^2} - {\mu(x)\over 2
h_{xm}} \right) e^{i h_{xm} a(x) p } \delta_{x', x-h_{xm}}
-{\sigma(x)^2\over h_{xm}^2} \delta_{x' x} - i p {b(x)}\delta_{x' x} \notag \\
& = {\sigma(x, p, h_{xm})^2\over 2} \Delta^m_x(x, x') + \mu(x, p,
h_{xm}) \nabla^m_x(x, x') + \kappa(x,
p, h_{xm}) \delta_{x x'} \notag \\
\label{eq_fourierl}
\end{align}
where
\begin{align}
&\sigma(x, p, h_{xm})^2 = \sigma(x)^2 \cos (h_{xm} a(x) p) - i  {\sin( h_{xm} p a(x))} \mu(x) \label{eq_sigma2}\\
&\mu(x, p, h_{xm}) = \mu(x) \cos (h_{xm} a(x) p) - i  {\sin( h_{xm} p a(x)) \over h_{xm}} \sigma(x)^2, \label{eq_mu2} \\
&\kappa(x, p, h_{xm}) = - i {\sin( h_{xm} p a(x)) \over h_{xm}}
\mu(x) + \sigma(x)^2 (\cos (h_{xm} a(x) p)-1) - i p b(x) . \label{eq_kappa2}
\end{align}

Notice that these three functions admit continuations as entire
analytic functions in $p$. Let us denote them with
$\sigma(x, z, h_{xm})^2, \mu(x, z, h_{xm}), \kappa(x, z, h_{xm})$.
Let $\hat\LLL^{m}(x, x'; z)$ be the analytic continuation of the
operator above and consider the Fourier transformed joint
propagator defined as follows:
\begin{align}
\hat u_m(x, x'; z; t) =& {1\over h_{xm}} \sum_{q=1}^\infty
\sum_{\gamma\in \Gamma_m : \gamma_0 = x, \gamma_{q} = x'}
 \int_0^{t} ds_1 \int_{s_1}^{t} ds_2 ...
\int_{s_{q-1}}^{t} d s_{q} \rho(s_1, s_2, ..., s_{q-1} ; z; \gamma)
\notag \\ \label{eq_upathint}
\end{align}
where
\begin{align}
\rho(s_1, s_2, ..., s_{q-1} ; z; \gamma)=& {1\over h_{xm}} e^{
-s_{1} \hat \LLL^m (\gamma_0, \gamma_0; z)} \notag
\prod_{j=1}^{q-1} \hat \LLL^m(\gamma_{j-1}, \gamma_j; z)
e^{-(s_{j+1}- s_{j}) \hat \LLL^m (\gamma_j, \gamma_j; z)} ds_1 ...
ds_{q}. \notag\\ \label{eq_upathint}
\end{align}
and $s_{q} = t$.

If $K>0$, let us introduce the constants
\begin{equation}
\Sigma_0 = \inf_{x \in X_m} \sigma(x, 0, h_{xm}), \;\;\;\;\;
M(K) = \sup_{\begin{matrix} x \in X_m \\ \lvert z \lvert < K
\end{matrix}} \lvert \mu (x, z, h_{xm}) \lvert
\end{equation}
and
\begin{equation}
\Sigma_1(K) = \sup_{\begin{matrix} x \in X_m \\ \lvert z \lvert < K
\end{matrix}} \sqrt{ \lvert\sigma(x, z, h_{xm})\lvert^2
+ \lvert\sigma(x, z, h_{xm}) - \sigma(x, 0, h_{xm}) \lvert^2 + h_{xm}
\lvert \mu (x, z, h_{xm})\lvert }.
\end{equation}
Due to our assumptions, $\Sigma_0>0$ and $\Sigma_1(K), M(K)
<\infty$.

A {\it symbolic path}  $\gamma = \{\gamma_0, \gamma_1, \gamma_2,
.... \}$ is an infinite sequence of sites in $X_m$ such that
$\gamma_j \neq \gamma_{j-1}$ for all $j=1, ...$. Let $\Gamma_m$ be
the set of all symbolic paths in $X_m$. The kernel of the diffusion
process admits the following representation in terms of a summation
over symbolic paths:
\begin{align}
\hat u_m(x, y; z; t) =& \sum_{q=1}^\infty 2 ^{-q}
\sum_{
\begin{matrix}
\gamma\in \Gamma_m : \gamma_0 = x, \gamma_q = y \\
\lvert\gamma_j - \gamma_{j-1}\lvert = 1 \;\;\forall j\ge1
\end{matrix}}
 W_m(\gamma, q, z, t)
 \label{eq_udiff1}
 \end{align}
where
\begin{align}
W_m(&\gamma, q, z, t) = {1\over h_{xm}} \int_0^{t} ds_1
\int_{s_1}^{t} ds_2 ... \int_{s_{q-1}}^{t} d s_{q} e^{ (t - s_q)
\hat {\mathcal L}^m(\gamma_q, \gamma_q; z)} \prod_{j=0}^{q-1}  2 \bigg(
\hat  {\mathcal L}^m(\gamma_j, \gamma_{j+1}; z) e^{
(s_{j+1} - s_{j}) \hat {\mathcal L}^m(\gamma_j, \gamma_{j}; z)}
 \bigg) \label{eq_wdiff}
\end{align}
with $s_0 = 0$.

\begin{figure}
\begin{center}
    \includegraphics[width = 12cm]{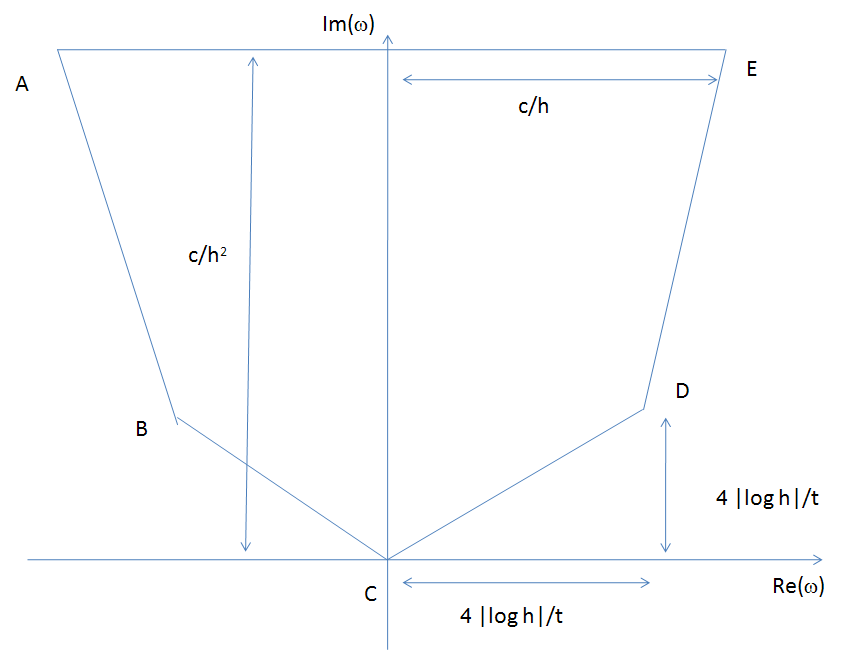}
    \caption{Contour of integration for the integral in (\ref{eq_greenfunc}).
    ${\mathcal C}_+$ is the countour joining the point $D$ to the
    points $E, A, B$. ${\mathcal C}_-$ is
    the countour joining the point $B$ to $C$ to $D$.}
    \label{fig_contour2}
\end{center}
\end{figure}

Let us introduce the following Green's function:
\begin{equation}
G_{m}(x, y; z, \omega) = \int_0^\infty \hat u_{m}(x, y; z, t)
e^{-i\omega t} dt = h_{xm}^{-1} {1\over \hat {\mathcal L}^m(z) +
i\omega}(x, y).
\end{equation}
The Fourier transformed propagator can be re-obtained from the
Green's function by evaluating the following contour integral:
\begin{equation}
\hat u_{m}(x, y; z, t) = \int_{{\mathcal C}_-} {d\omega \over 2 \pi}
G_{m}(x, y; z, \omega) e^{i\omega t} + \int_{{\mathcal C}_+}
{d\omega \over 2 \pi} G_{m}(x, y; z, \omega) e^{i\omega t}.
\label{eq_greenfunc}
\end{equation}
Here, ${\mathcal C}_+$ is the contour joining the point $D$ to the
points $E, A, B$ in Fig. \ref{fig_contour2}, while ${\mathcal C}_-$
is the contour joining the point $B$ to $C$ to $D$. The intent of the
design of the contour of integration is to ensure that each
point $\omega$ on the upper path ${\mathcal C}_+$ is separated from
the spectrum of $\mathcal L$.

\begin{lemma} \label{lemma_cplus} For all $K>0$, there are an integer
$m(K)$ and a constant $c(K)$ such that
\begin{equation}
\bigg\lvert\int_{{\mathcal C}_+} {d\omega \over 2 \pi} G_{m}(x, y;
z,  \omega) e^{i\omega t}\bigg\lvert \leq c(K) h^2 .
\end{equation}
for all $z$ such that $\lvert z \lvert <K$.
\end{lemma}
\begin{proof}
The proof is based on the geometric series expansion
\begin{align}
G_{m}(z, \omega) &= h_m^{-1} {1\over {\mathcal L}^m(z) + i\omega} =
h_m^{-1} \sum_{j=0}^\infty {1\over {1\over2}\sigma(0,
h_{xm})^2\Delta^m + i\omega} \notag \\ & \hskip2cm \bigg[ \bigg(
{1\over2}(\sigma(z, h_{xm})^2 - \sigma(0, h_{xm})^2)\Delta^m +
\mu(z, h_{xm})\nabla^m  + \kappa(z, h_{xm}) \bigg)
{1\over {1\over2}\sigma(0, h_{xm})^2\Delta^m + i\omega} \bigg]^j \notag \\
\label{geoseries}
\end{align}
Here $\sigma(z, h_{xm})^2$ and $\mu(z, h_{xm})$ are the
multiplication operators by $\sigma(x, z, h_{xm})^2$ and $\mu(x, z,
h_{xm})$, respectively. Convergence for $\omega\in{\mathcal C}_+$
can be established by means of a Kato-Rellich type relative bound,
see \cite{Kato}. More precisely, for any $\alpha>0$, one can find a
$\beta>0$ such that the operators $\nabla^m$ and $\Delta^m$ satisfy
the following relative bound estimate:
\begin{equation}
\lvert\lvert \nabla^m f\lvert\lvert_2 \leq \alpha \lvert\lvert
\Delta^m f \lvert\lvert_2 + \beta \lvert\lvert f \lvert\lvert_2.
\end{equation}
for all periodic functions $f$ and all $m\ge m_0$. This bound can be
derived by observing that $\nabla^m$ and $\Delta^m$ can be
diagonalized simultaneously by a Fourier transform and by observing
that for any $\alpha>0$, one can find a $\beta>0$ such that
\begin{equation}
\bigg\lvert {\sin h_{xm} k \over h_{xm}}  \bigg\lvert \leq \alpha
\bigg\lvert {\cos h_{xm} k - 1 \over h_{xm}^2} \bigg\lvert + \beta
\end{equation}
for all $m\ge m_0$ and all $k\in B_m$.

Under the same conditions, we also have that
\begin{align} &\bigg\lvert\bigg\lvert
\bigg( {1\over2}(\sigma(z, h_{xm})^2 - \sigma(0, h_{xm})^2)\Delta^m
+ \mu(z, h_{xm})\nabla^m \bigg) f\bigg\lvert\bigg\lvert_2  \notag
\\ &\hskip5cm \leq {2 \alpha \over \Sigma_0^2} \big( M(K) + 2 \Sigma_1(K)^2\big)
\bigg\lvert\bigg\lvert {1\over2} \sigma(0)^2 \Delta^m f
\bigg\lvert\bigg\lvert_2 + \beta \lvert\lvert f \lvert\lvert_2.
\end{align}
Hence
\begin{align} \bigg\lvert\bigg\lvert
\mu & \nabla^m {1\over {1\over2}\sigma^2\Delta^m + i\omega}
f\bigg\lvert\bigg\lvert_2 \notag \\ & \leq \bigg(
{1\over2}(\sigma(z, h_{xm})^2 - \sigma(0, h_{xm})^2)\Delta^m +
\mu(z, h_{xm})\nabla^m \bigg) \bigg\lvert\bigg\lvert {1\over2}
\sigma^2 \Delta^m {1\over {1\over2}\sigma^2\Delta^m + i\omega} f
\bigg\lvert\bigg\lvert_2 \notag \\ & \hskip9cm + \beta
\bigg\lvert\bigg\lvert {1\over {1\over2}\sigma^2\Delta^m + i\omega}
f \bigg\lvert\bigg\lvert_2 < 1
\notag \\
\end{align}
where the last inequality holds if $\omega\in{\mathcal C}_+$, if
$\alpha$ is chosen sufficiently small and if $m$ is large enough. In
this case, the geometric series expansion in
(\ref{geoseries}) converges in $L^2$ operator norm. The uniform norm
of the kernel $\lvert G_{m}(x, y; \omega) \lvert$ is pointwise
bounded from above by $h_m^{-2}$.

Since the points $B$ and $D$ have imaginary part of height $
4{\lvert \log h_m \lvert \over t}$, the integral over the contour
${\mathcal C}_+$ converges also and is bounded from above by $c
h_m^2$ in uniform norm.

\end{proof}

\begin{lemma} For all $K>0$, if $q\ge {e^2
 t\over 2 h_m^2}(2 \Sigma_1(K)^2 + M(K))$ we have that
\begin{equation}
W_m(\gamma, q; z, t) \leq \sqrt{q\over 2\pi} \exp\left(-{\Sigma_0^2
t\over 2} - q\right). \label{eq_wbound}
\end{equation}
\end{lemma}
\begin{proof}
Let us define the function
\begin{equation}
\phi(t) = {\Sigma_1(K)^2 \over 2 h_m^2} \; e^{-{\Sigma_0^2 t \over 2
h_m^2}} \; 1(t\ge 0)
\end{equation}
where $1(t\ge0)$ is the characteristic function of $\RRR_+$. We have
that
\begin{equation}
W_m(\gamma, q; z, t) \leq \phi^{\star q}(t)
\end{equation}
for all $z$ such that $\lvert z\lvert < K$, where $\phi^{\star q}$
is the $q-$th convolution power, i.e. the $q-$fold convolution
product of the function $\phi$ by itself. The Fourier transform of
$\phi(t)$ is given by
\begin{equation}
\hat\phi(\omega) = {\Sigma_1(K)^2 \over 2 h_m^2} \int_0^\infty
e^{-i\omega t - {\Sigma_0^2 t\over 2 h_m^2}} dt = {\Sigma_1(K)^2
\over 2 i \omega h_m^2 + \Sigma_0^2}.
\end{equation}
The convolution power is given by the following inverse Fourier
transform:
\begin{equation}
\phi^{\star q}(t) = \int_{-\infty}^\infty \hat \phi(\omega)^q e^{i \omega t}
{d\omega\over 2\pi} = \left( {\Sigma_1(K)\over\Sigma_0} \right)^{2 q}
\int_{-\infty}^\infty\left( 1 + {2 i \omega h_m^2 \over \Sigma_0^2}
\right)^{-q} e^{i \omega t} {d\omega\over 2\pi}.
\end{equation}
Introducing the new variable $z = 1 + {2i\omega h_m^2\over
\Sigma_0^2}$, the integral can be recast as follows
\begin{equation}
\phi^{\star q}(t) = {\Sigma_0^{2-2q} \Sigma_1(K)^{2q}\over 4\pi i
h_m^2} \lim_{R\to\infty} \int_{{\mathcal C}_R} z^{-q} \exp\left(
{\Sigma_0^2 t\over 2 h_m^2} (z-1) \right) dz \label{eq_intc}
\end{equation}
where ${\mathcal C}_R$ is the contour in Fig. \ref{fig_contour1}.
Using the residue theorem and noticing that the only pole of the
integrand is at $z = 0$, we find
\begin{equation}
\phi^{\star q}(t) = {1\over (q-1)!} \left({\Sigma_1(K)^2 t \over 2
h_m^2}\right)^q \exp\left( - \Sigma_0^2 t \over 2 h_m^2 \right).
\end{equation}
Making use of Stirling's formula $q! \approx \sqrt{2\pi}
q^{q+{1\over2}} e^{-q}$, we find
\begin{equation}
\phi^{\star q}(t) \approx \sqrt{q\over 2\pi} \exp\left( -
{\Sigma_0^2 t \over 2 h_m^2} + q\log{\Sigma_1(K)^2 t \over 2 h_m^2}
+ q (1-\log q) \right).
\end{equation}
If $\log q \ge \log {\Sigma_1(K)^2 t \over 2 h_m^2} +2$, then we
arrive at the bound in (\ref{eq_wbound}).

\begin{figure}
\begin{center}
    \includegraphics[width = 12cm]{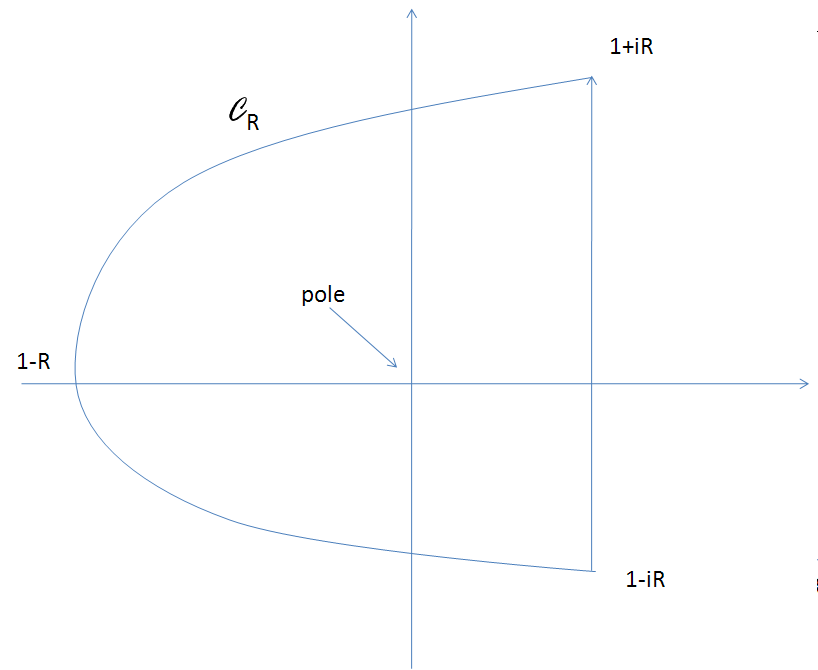}
    \caption{Contour of integration ${\mathcal C}_R$ for the integral in (\ref{eq_intc}).}
    \label{fig_contour1}
\end{center}
\end{figure}

\end{proof}

Let us fix a $K>0$. To prove the theorem, it suffices to
consider the case $m' = m+1$ for all values of $m$ above $m_0$. In
fact, given this particular case, the general statement can be
derived with an iterative argument. To this end, we introduce a
renormalization group transformation based on the notion of
decorating path.

\begin{definition}{\bf (Decorating Paths.)}
 Let $m\ge m_0$ and let $\gamma = \{y_0, y_1, y_2,
.... \}$ be a symbolic sequence in $\Gamma_m$. A {\it decorating
path around $\gamma$} is defined as a symbolic sequence $\gamma' =
\{y_0, y_1', y_2', .... \}$ with $y_i' \in h_{m+1} \ZZZ$ containing
the sequence $\gamma$ as a subset and such that if $y_j' = y_i$ and
$y_k' = y_{i+1}$, then all elements $y_n'$ with $j < n < k$ are such
that $\lvert y_n' - y_j' \lvert \leq h_{m+1}$. Let ${\mathcal
D}_{m+1}(\gamma)$ be the set of all decorating sequences around
$\gamma$. The decorated weights are defined as follows:
\begin{equation}
{\tilde W}_{m}(\gamma, q; z;  t) = \sum_{q' = q}^\infty \sum_{
\begin{matrix}
\gamma' \in {\mathcal D}_{m+1}(\gamma)\\
\gamma'_{q'} = \gamma_q
\end{matrix}
} W_{m+1}(\gamma', q';  z; t).
\end{equation}
Let us notice that these weights are not positive or even real
values unless $z=0$, but rather the depend analytically on the
variable $z$. Finally, let us introduce also the following Fourier
transform:
\begin{equation}
\hat W_{m}(\gamma, q;  z; \omega) = \int_0^\infty W_{m}(\gamma, q;
 z; t) e^{i\omega t} dt, \;\;\;\; \hat {\tilde W}_{m}(\gamma, q;
 z; \omega) = \int_0^\infty \tilde W_{m}(\gamma, q;  z; t) e^{i\omega t}
dt.
\end{equation}
\end{definition}

\begin{notation}
In the following, we set $h = h_{x, m+1}$ so that $h_{xm} = 2 h$. We
also use the Landau notation $O(h^n)$ to indicate a function $f(h)$
such that $ h^{-n} f(h) $ is bounded in a neighborhood of $0$.
\end{notation}

\begin{lemma}  Let $x, y\in X_m$ and let ${\mathcal C_-}$ be an
integration contour as in Fig. \ref{fig_contour2}. Then
\begin{equation}
\bigg\lvert \bigg( \int_{{\mathcal C}_-}  2 G_{m+1}(x, y;  z;
\omega) - G_m(x, y;  z; \omega) \bigg) e^{i\omega t} {d\omega\over
2\pi} \bigg\lvert = O(h^2).
\end{equation}
\end{lemma}

\begin{proof}

We have that
\begin{align}
2 G_{m+1}(x, y;  z; \omega) - G_m(x, y;  z; \omega) = & {1\over h}
\sum_{q=1}^\infty 2 ^{-q} \sum_{
\begin{matrix}
\gamma\in \Gamma_m : \gamma_0 = x, \gamma_{q} = y\\
\lvert\gamma_j - \gamma_{j-1}\lvert = 1 \forall j\ge1
\end{matrix}}
\left( 2 \hat {\tilde W}_m (\gamma, q;  z; \omega) - \hat W_m
(\gamma, q;  z; \omega)\right). \label{eq_udiff}
\end{align}
The number of paths over which the summation is extended is
\begin{equation}
N(\gamma, q; x, y) \equiv \sharp\{\gamma\in \Gamma_m : \gamma_0 = x,
\gamma_{q} = y, \lvert\gamma_j - \gamma_{j-1}\lvert = 1 \forall
j\ge1 \} = \left(\begin{matrix} q \\ {q\over2} + k
\end{matrix}\right)
\end{equation}
where $k = {\lvert y - x \lvert \over h_{xm}}.$ Applying Stirling's
formula we find
\begin{equation}
N_\gamma \lesssim 2^q \sqrt{2\over \pi q}.
\end{equation}
Hence
\begin{align}
&\bigg\lvert \int_{{\mathcal C}_-} \bigg( 2 G_{m+1}(x, y;  z;
\omega) - G_m(x, y;  z; \omega) \bigg) e^{i\omega t}  {d\omega\over
2\pi}
\bigg \lvert \notag \\
&\leq {c\over h}
 \sum_{q=1}^\infty \sqrt{1\over q} \max_{
\begin{matrix}
\gamma\in \Gamma_m : \gamma_0 = x, \gamma_{q} = y\\
\lvert\gamma_j - \gamma_{j-1}\lvert = 1 \forall j\ge1
\end{matrix}}
\bigg \lvert \int_{{\mathcal C}_-} \bigg(2 \hat {\tilde W}_m
(\gamma, q; z, \omega) - \hat W_m (\gamma, q;  z; \omega)\bigg)
e^{i\omega t}  {d\omega\over 2\pi} \bigg\lvert. \notag
\\\label{eq_udiff}
\end{align}
for some constant $c\approx\sqrt{2\over\pi} >0$. It suffices to
extend the summation over $q$ only up to
\begin{equation}
q_{\rm max} \equiv {e^2 \Sigma_1(K)^2 t\over 2 h^2}.
\end{equation}
To resum beyond this threshold, one can use the previous lemma. More
precisely, we have that
\begin{align}
&\bigg\lvert \int_{{\mathcal C}_-} \bigg( 2 G_{m+1}(x, y;  z;
\omega) - G_m(x,
y;  z; \omega) \bigg) e^{i\omega t}  {d\omega\over 2\pi} \bigg\lvert \notag \\
&\leq {c \sqrt{q_{\rm max}}\over h} \max_{
\begin{matrix}
q, \gamma\in \Gamma_m : \gamma_0 = x, \gamma_{q} = y\\
\lvert\gamma_j - \gamma_{j-1}\lvert = 1 \forall j\ge1
\end{matrix}}
\bigg \lvert \int_{{\mathcal C}_-} \bigg( 2 \hat {\tilde W}_m
(\gamma, q; z, \omega) - \hat W_m (\gamma, q; z, \omega) \bigg)
e^{i\omega t}  {d\omega\over 2\pi}  \bigg\lvert. \notag
\\\label{eq_udiff2}
\end{align}

Let us introduced the following abbreviated notations:
\begin{align}
& v(x, z) = \sigma(x,  z;, h)^2 \\
& m(x, z) = \mu(x,  z;, h) \\
& k(x, z) = \kappa(x,  z;, h)
\end{align}
where for the sake of keeping formulas short, we omit to denote the
$ h = h_{x, m+1}$ dependencies. To evaluate the resummed weight
function, let us form the matrix
\begin{equation} \bar {\mathcal L}(x, z) =
\left(
\begin{matrix}
-{v\left(x+h, z\right)\over h^2} + \kappa\left(x+h, z\right) && {v\left(x+ h, z\right)\over 2
h^2} -
{m(x+h, z)\over 2 h} && 0 \\
{v\left(x, z\right)\over 2 h^2}  +{m(x, z)\over 2h} && -{ v\left(x , z\right)\over h^2} +\kappa\left(x, z\right)  && {v\left(x, z\right)
\over 2 h^2} - {m(x, z) \over 2 h} \\
0 && {v\left(x - h, z\right)\over 2 h^2} + {m(x-h, z)\over 2 h} &&
-{v\left(x - h , z\right) \over h^2} \kappa\left(x-h, z\right)
\end{matrix} \right)
\end{equation}
and decompose it as follows:
\begin{equation}
\bar {\mathcal L}(x, z) = {1\over h^2} \bar {\mathcal L}_0(x, z) +
{1\over h} \bar {\mathcal L}_1(x, z) + \bar {\mathcal L}_2(x, z) + h
\bar {\mathcal L}_3(x, z).
\end{equation}
where
\begin{equation} \bar {\mathcal
L}_0(x, z) = \left(
\begin{matrix}
-v(x, z) && {1\over2} {v(x, z)} && 0 \\
{1\over2} {v(x, z)} && -{v(x, z)} && {1\over2} {v(x, z)} \\
0 && {1\over2} {v(x, z)} && -{v(x, z)}
\end{matrix}
\right),
\end{equation}
\begin{equation}
\bar {\mathcal L}_1(x, z) = \left(
\begin{matrix}
-\nabla_x^m v(x, z) && {1 \over2} \nabla_x^m v(x, z) - {1\over2} m(x, z) && 0 \\
{1\over2} m(x, z) && 0 && -{1\over2} m(x, z) \\
0 && - {1\over2} \nabla_x^m v(x, z) + {1\over2}  m(x, z) &&
\nabla_x^m v(x, z)
\end{matrix}\right),
\end{equation}
\begin{equation}
\bar {\mathcal L}_2(x, z) = \left(
\begin{matrix}
-{1\over2} \Delta_x^m v(x, z) + k(x, z) && {1\over4} \Delta_x^m v(x, z) - {1\over2} \nabla_x^m m(x, z)&& 0 \\
0 && k(x, z) && 0 \\
0 && {1\over4} \Delta_x^m v(x, z) - {1\over2} \nabla_x^m m(x, z) &&
-{1\over2} \Delta_x^m v(x, z) + k(x, z)
\end{matrix}\right).
\end{equation}
and
\begin{equation}
\bar {\mathcal L}_3(x, z) = \left(
\begin{matrix}
\nabla^m_x k(x, z)  && - {1\over4} \Delta_x^m m(x, z)&& 0 \\
0 && 0 && 0 \\
0 && {1\over4} \Delta_x^m m(x, z) && - \nabla^m_x k(x, z)
\end{matrix}\right).
\end{equation}
Recall that all functions and operators above depend also on $h =
h_{x, m+1}$.

Let us introduce the sign variable $\tau = \pm1$, the functions
\begin{align}
\phi_{0}(t, x, z, \tau) &\equiv 2 \bar {\mathcal L}^m_z(x, x+2\tau
h; z) e^{ t \hat {\mathcal
L}_z^{m}(x, x; z)} 1(t\ge0) \\
\phi_{1}(t, x, z, \tau) &\equiv 2 \bar {\mathcal L}^{m+1}_z(x+\tau
h, x+2\tau h; z) e^{t \bar {\mathcal L}(x, z)}(x, x+\tau h) 1(t\ge0)
\end{align}
and their Fourier transforms
\begin{align}
\hat\phi_{0}(\omega, x, z, \tau) &= \left({v(x, z)\over 4 h^2} +
\tau { m(x, z)\over 2 h}\right) \left( {v(x, z)\over 4 h^2} - k(x,
z) + i \omega \right)^{-1} \notag \\ \hat\phi_{1}(\omega, x, z,
\tau) &= \left({v(x, z)\over h^2} + \tau  {m(x, z) + \nabla_x^m v(x,
z)\over h} + {\Delta^m_x v(x, z)+ \nabla^m_x  m(x, z)\over2} +
{\Delta^m_x
 m(x, z)\over2} \tau h + O(h^2)\right) \notag \\
&\hskip8cm < x \lvert \left( -\bar {\mathcal L}(x, z) + i\omega
\right)^{-1}\lvert x+\tau h >.
\end{align}
where
\begin{equation}
\lvert x> = \left(\begin{matrix} 0 \\ 1 \\
0 \end{matrix}\right), \;\;\;\;{\rm and}\;\;\;\; \lvert x+\tau h> =
\left( \begin{matrix} \delta_{\tau, 1} \\ 0 \\
\delta_{\tau, -1} \end{matrix}\right).
\end{equation}
We also require the functions
\begin{equation}
\psi_{0}(t, x, z) \equiv e^{ t {\mathcal L}^m_z(x, x; z) } 1(t\ge0),
\;\;\;\;\;\; \psi_{1}(t, x) \equiv e^{t \bar {\mathcal L}(y; h)}(x,
x) 1(t\ge0)
\end{equation}
and the corresponding Fourier transforms
\begin{align}
\hat\psi_{0}(\omega, x, z) = \left( {v(x, z)\over 4 h^2}  + i \omega
\right)^{-1}, \;\;\;\;\; \hat\psi_{1}(\omega, x, z) =  <x \lvert
\left( -\bar {\mathcal L}(x, z) + i\omega \right)^{-1}\lvert x>.
\end{align}

If $\gamma$ is a symbolic sequence, then
\begin{align}
\hat W_m(\gamma, q; z, \omega) &= \hat\psi_{0}(\omega,\gamma_q, z)
\prod_{j=0}^{q-1} \hat\phi_0(\omega;
\gamma_j, z, {\rm sgn}(\gamma_{j+1}-\gamma_j))\\
 \hat {\tilde W}_m(\gamma, q; z; \omega)
&= \hat\psi_{1}(\omega, \gamma_q, z) \prod_{j=0}^{q-1}
\hat\phi_1(\omega; \gamma_j, p, {\rm sgn}(\gamma_{j+1}-\gamma_j)).
\end{align}

Let us estimate the difference between the functions
$\hat\phi_{1}(\omega, x, z, \tau)$ and $\hat\phi_{2}(\omega, x, z,
\tau)$ assuming that $\omega$ is in the contour ${\mathcal C}_-$ in
Fig. \ref{fig_contour1}. Retaining only terms up to order up to
$O(h^3)$, we find
\begin{align}
\hat\phi_{0}(\omega, x, z, \tau) =& 1 + {2m(x, z)\tau h \over v(x,
z)} + {4 h^2 \over v(x, z)}(k(x, z) - i \omega) + \notag \\
&\hskip2cm{8 m(x, z) \tau h^3 \over v(x, z)^2} (k(x, z) - i \omega)
+ {16 h^4 \over v(x, z)^2}(k(x, z) - i \omega)^2
+ O(h^5). \notag \\
\end{align}
A lengthy but straightforward calculation which is best carried out
using a symbolic manipulation program, gives
\begin{align}
&\hat\phi_{1}(\omega, x, z, \tau) = 1+{2m(x, z) \tau h \over v(x,
z)}+{ 4
 h^2 \over v(x, z)}(k(x, z) - i \omega)  - \big[8 m(x, z)  - \nabla_x^m v(x, z) \big]
{i\omega\tau h^3 \over v(x, z)^2} \notag \\
&\hskip6cm+ R(x, z)\cdot h^3 \tau
+ h^4 P_0(x, z) + i\omega h^4 P_1(x, z) - {14\omega^2 h^4 \over v(x, z)^2} + O(h^5) \notag \\
\end{align}
where
\begin{align}
&R = {1\over 2 v^3} \bigg[ \big(\Delta^m_x m + 2 \nabla^m_x k\big)
v^2 + \big(16 m - 2 \nabla_x^m v \big) k v\notag
\\ &\hskip2cm - 4m^3  + 2 \nabla_x^m v m^2 - 2
v \nabla_x^m v \nabla_x^m m - \big(2m \nabla_x^m m + v \Delta^m_x v
-
2\big(\nabla_x^m v\big)^2\big) m \bigg] . \notag \\
&P_0 = {1 \over v^3} \bigg[\bigg((2 m - 2\nabla_x^m v) \nabla_x^m k +
14 k^2 + (-4 \nabla_x^m m - \Delta_x^m v) k\bigg) v + \bigg(-4 m^2+2
m \nabla_x^m v + 2 (\nabla_x v)^2 \bigg) k \bigg]
\notag \\
&P_1 = {1\over v^3} \bigg[ 28 k - 4m^2 + 2 m \nabla_x^m v -4 v
\nabla_x m - v \Delta^m_x v
+ 2 \big(\nabla^m_x v\big)^2  \bigg].  \notag \\
\end{align}
For simplicity, we are not denoting here the dependency of all
functions on $(x, z)$. We have that
\begin{align}
& \sum_{j=0}^{q-1} \bigg(\log \hat\phi_0(\omega; \gamma_j, z, {\rm
sgn}(\gamma_{j+1}-\gamma_j)) -
 \log \hat\phi_1(\omega; \gamma_j, z, {\rm
sgn}(\gamma_{j+1}-\gamma_j))\bigg)
\notag \\
&= \sum_{j=0}^{q-1} \bigg(  {i \omega \nabla_x^m v(\gamma_j, z)\over
v(\gamma_j, z)^2} + R(\gamma_j, z) \bigg) h^3 {\rm
sgn}(\gamma_{j+1}-\gamma_j) + \big(\lvert\lvert P_0\lvert\lvert_{K,
\infty} + \lvert \omega\lvert \lvert\lvert P_1\lvert\lvert_{K,
\infty} + 2 \lvert \omega\lvert^2 \lvert\lvert v^{-2}
\lvert\lvert_{K,
\infty} \big) O(h^4 q) \notag \\
&= i\omega h^2 \bigg({1\over v(\gamma_0,
z)} - {1\over v(\gamma_{q}, z)} \bigg) + h^2 \big(G(\gamma_q, z) - G(\gamma_0, z)\big) +
\big(\lvert\lvert P_0\lvert\lvert_{K, \infty} + \lvert \omega\lvert
\lvert\lvert P_1\lvert\lvert_{K, \infty} + 2 \lvert \omega\lvert^2
\lvert\lvert v^{-2} \lvert\lvert_{K,
\infty} \big) O(h^4 q) \notag \\
\end{align}
where $G(x, z)$ is a primitive of $R(x, z)$, i.e.
\begin{equation}
G(x, z) = \int^x_{-L_x} R(x', z) dx'.
\end{equation}

We conclude that there is a constant $c(K)>0$ such that
\begin{equation}
\bigg\lvert \int_{{\mathcal C}_-} \bigg( \prod_{j=0}^{q-1}
\hat\phi_0(\omega; \gamma_j, z, {\rm sgn}(\gamma_{j+1}-\gamma_j)) -
\prod_{j=0}^{q-1} \hat\phi_1(\omega; \gamma_j, z, {\rm
sgn}(\gamma_{j+1}-\gamma_j)) \bigg) e^{i\omega t}  {d\omega\over
2\pi} \bigg\lvert \leq c(K) h^\gamma.
\end{equation}
for all $q\leq q_{\max}$. Here we use the decay of $e^{i\omega t}$
in the upper half of the complex $\omega$ plane to offset the
$\omega$ dependencies in the integrand. Similar calculations lead to
the following expansions:
\begin{equation}
\hat\psi_{0}(\omega, x, z) = {4 h^2 \over v(x, z)}  + O(\omega h^4),
\;\;\;\;\; \hat\psi_{1}(\omega, x, z) =  {2 h^2 \over v(x, z)} +
O(\omega h^4) = {1\over 2} \hat\psi_{0}(\omega, x, z) + O(\omega
h^4).
\end{equation}

Since $q<c h^{-2}$ and $\omega \leq \lvert \log h \lvert$, we find
\begin{align}
\bigg\lvert \int_{{\mathcal C}_-} \bigg(  2 G_{m+1}(x, y; z, \omega)
- G_m(x, y; z, \omega) \bigg) e^{i\omega t}  {d\omega\over 2\pi}
\bigg\lvert \leq c(K) { \sqrt{q_{\rm max}} \over h }  h^{2+\gamma} \leq c
h^\gamma. \label{eq_udiff3}
\end{align}

By differentiating with respect to time in equation
(\ref{eq_greenfunc}), we find that
\begin{equation}
{\partial \over \partial t} u_{m}(x, y; t) = \int_{{\mathcal C}_-}
 i\omega G_{m}(x, y;  z; \omega) e^{i\omega t} {d\omega \over 2 \pi}
+ \int_{{\mathcal C}_+}  i\omega  G_{m}(x, y;
 z; \omega) e^{i\omega t} {d\omega \over 2 \pi}. \label{eq_greenfunc}
\end{equation}
All the derivations above carry through. We conclude that
for some constant $c_1(K)>0$ and all $\lvert z \lvert < K$ we have that
\begin{equation}
\bigg\lvert\int_{{\mathcal C}_+} {d\omega \over 2 \pi} i\omega
G_{m}(x, y;  z; \omega) e^{i\omega t}\bigg\lvert \leq c_1(K) h^\gamma
\end{equation}
and also
\begin{align}
\bigg\lvert \int_{{\mathcal C}_-} i\omega \bigg(  2 G_{m+1}(x, y;
 z; \omega) - G_m(x, y;  z; \omega) \bigg) e^{i\omega t}  {d\omega\over
2\pi} \bigg\lvert \leq c { \sqrt{q_{\rm max}} \over h }  h^4 \leq
c_1(K) h^\gamma. \label{eq_udiff3}
\end{align}
Hence, the first time derivatives of the kernel satisfy a similar
Cauchy convergence condition as the kernel itself.

Let us notice that, if $m$ is large enough, then
the operator $\hat {\mathcal L}^m(z)$ is relatively bounded with respect to
${\mathcal L}^m$ in the following sense:
\begin{equation}
\lvert\lvert {\mathcal L} f \lvert\lvert
\leq
\alpha(K) \lvert\lvert \hat {\mathcal L}(z) f \lvert\lvert
+ \beta(K) \lvert\lvert f \lvert\lvert
\end{equation}
for some $\alpha(K), \beta(K)>0$ and all $\lvert z \lvert <K$.
By inspection of the terms in $\hat {\mathcal L}(z)$
one sees that such bound holds if $h_{xm}$ is small enough, i.e. for
$m$ larger than a threshold depending on $K$. A similar relative bound
holds for the adjunct operator ${\mathcal L}^*$ with respect to $\hat {\mathcal L}(z)^*$.
This concludes the proof of Theorem \ref{kernel_convergence}.

\end{proof}

\section{Explicit Euler Scheme \label{SecEuler}}

\begin{figure}
\begin{center}
    \includegraphics[width = 12cm]{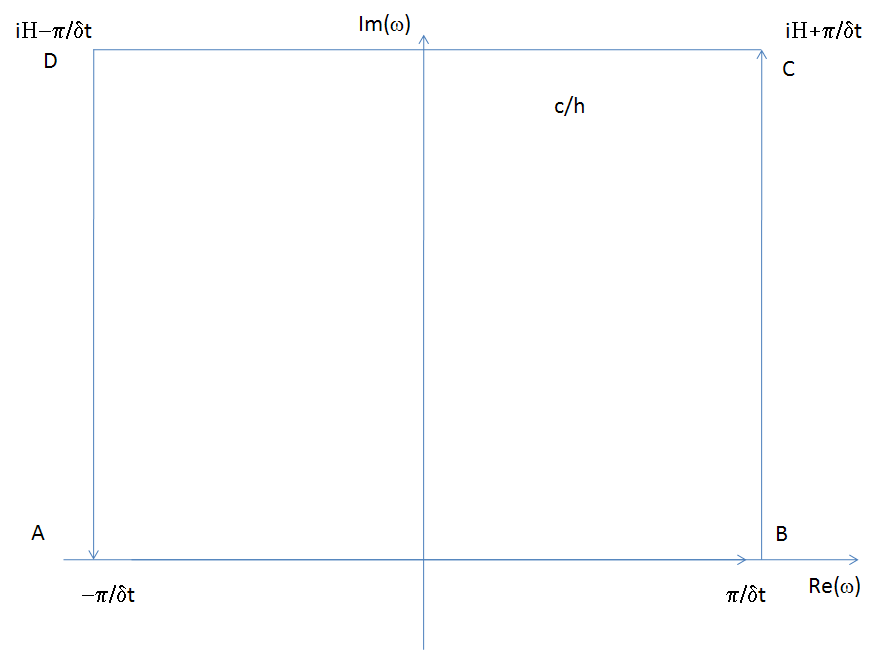}
    \caption{Contour of integration for the integral in (\ref{eq_greenfuncdt}).}
    \label{fig_contour3}
\end{center}
\end{figure}

In this section we prove Theorem \ref{kernel_euler}.
A path-wise expansion for the time-discretization of the Fouier transformed kernel
has the form
\begin{align}
\hat u_m^{\delta t}(x, y; z, t) =& {1\over h_{xm}} \sum_{q=1}^\infty
\sum_{\gamma\in \Gamma_m : \gamma_0 = x, \gamma_{q} = y}
 \sum_{k_1 = 1}^{N} \sum_{k_2 = k_1 + 1}^{N} ..
\sum_{k_{q} = k_{q-1} + 1}^{N}  \notag \\
& \bigg(1 + \delta t \hat \LLL_m(\gamma_0, \gamma_0; z)\bigg)^{k_{1}-1}
(\delta t)^q \prod_{j=1}^{q} \hat \LLL_m(\gamma_{j-1}, \gamma_j; z) \bigg(1
+ \delta t \hat \LLL_m(\gamma_j, \gamma_j; z)\bigg)^{k_{j+1} - k_{j} - 1}
\label{eq_upathint}
\end{align}
where $t_{q+1} = t$ and $k_{q+1}=N$. In this case, the propagator
can be expressed through a Fourier integral as follows:
\begin{equation}
\hat u_m^{\delta t}(x, y; z, t) = \int_{-{\pi\over\delta
t}}^{\pi\over\delta t} G_m^{\delta t}(x, y; z, \omega) e^{i\omega
t} {d\omega\over 2\pi}
\end{equation}
where
\begin{equation}
G_m^{\delta t}(x, y; z, \omega) = \delta t \sum_{j = 0}^{\infty}
\hat u_m^{\delta t}(x, y; z, j \delta t) e^{ - i\omega j
\delta t}.
\end{equation}
The Fourier transformed propagator can also be represented as the limit
\begin{equation}
\hat u_m^{\delta t}(x, y; z, t) = \lim_{H\to\infty} \int_{{\mathcal
C}_H} G_m^{\delta t}(x, y; z, \omega) e^{i\omega t} {d\omega\over
2\pi} \label{eq_greenfuncdt}
\end{equation}
where ${\mathcal C}_H$ is the contour in Fig. \ref{fig_contour3}.
This is due to the fact that the integral along the segments $BC$
and $DA$ are the negative of each other, while the integral over
$CD$ tends to zero exponentially fast as $\Im(\omega) \to \infty$,
where $\Im(\omega)$ is the imaginary part of $\omega$. Using
Cauchy's theorem, the contour in Fig. \ref{fig_contour3} can be
deformed into the contour in Fig. \ref{fig_contour2}. To estimate
the discrepancy between the time-discretized kernel and the
continuous time one, one can thus compare the Green's function along
such contour. Again, the only arc that requires detailed attention
is the arc $BCD$, as the integral over rest of the contour of
integration can be bounded from above as in the previous section.

Let $h = h_{xm}$ and let us introduce the two functions
\begin{align}
\phi_{0}(t, x, z, \tau) &\equiv 2 \hat {\mathcal L}_m(x, x+\tau h; z)
e^{ t \hat  {\mathcal L}_m(x, x; z)} 1(t\ge0), \\
\phi_{\delta t}(j, x, z, \tau) &\equiv 2 \hat {\mathcal L}_m(x, x+\tau h; z)
\big( 1 + \delta t \hat {\mathcal L}_m(x, x; z) \big)^{j-1}.
\end{align}
and the corresponding Fourier transforms
\begin{align}
\hat\phi_{0}(\omega, x, z, \tau) &= \int_0^\infty \phi_{0}(t, x, z,
\tau) e^{-i\omega t} {d\omega\over 2\pi} =  \left({v(x, z, h)\over h^2} +
\tau { m(x, z, h)\over h}\right) \left( {v(x, z, h)\over
h^2} + i \omega \right)^{-1} \\
\hat\phi_{\delta t}(\omega, x, z, \tau) &= \sum_{j=0}^{t\over \delta
t} \phi_{\delta t}(j, x, z, \tau) e^{-i\omega j\delta t} =
\left({v(x, z, h)\over h^2} + \tau {m(x, z, h)\over h}\right) \left(
e^{i\omega\delta t} - 1 + \delta t {v(x, z, h)\over h^2} \right)^{-1}.
\end{align}
We have that
\begin{align}
\hat\phi_{\delta t}(\omega, x, z, \tau) &= \left({v(x, z, h)\over h^2} +
\tau {m(x, z, h)\over h}\right) \left( i\omega  + {v(x, z, h)\over h^2} -
{\omega^2 \over 2} \delta t + O(\delta t^2)
\right)^{-1} \notag \\
& = \hat\phi_{0}(\omega, x, z, \tau) + {\omega^2\over 2 v(x, z, h)} h^2
\delta t + O(h^2 \delta t^2).
 = \hat\phi_{0}(\omega, x, z, \tau) +  O(h^4),
 \end{align}
where the last step uses the fact that $\delta t = O(h^2)$.

Let us also introduce the functions
\begin{align}
\psi_{0}(t, x, z, \tau) &\equiv e^{ t \hat {\mathcal L}_m(x, x; z)} 1(t\ge0),
\hskip1cm \psi_{\delta t}(j, x, \tau) \equiv \sum_{k = 1}^j \big( 1
+ \delta t \hat {\mathcal L}_m(x, x; z) \big)^{j-1}.
\end{align}
and the corresponding Fourier transforms
\begin{align}
\hat\psi_{0}(\omega, x, z, \tau) &= \left( {v(x, z, h)\over h^2} + i \omega
\right)^{-1}, \hskip1cm \hat\psi_{\delta t}(\omega, x, \tau) =
\left( e^{i\omega\delta t} - 1 + \delta t {v(x, z, h)\over h^2}
\right)^{-1}.
\end{align}
Again we find that
\begin{align}
\hat\psi_{0}(\omega, x, z, \tau) = \hat\psi_{\delta t}(\omega, x, z, \tau)
 +  O(h^4).
\end{align}

If $\gamma$ is a symbolic sequence, then let us set
\begin{align}
\hat W_m(\gamma, q; z, \omega) &= \hat\psi_{0}(\omega, \gamma_q)
\prod_{j=0}^{q-1} \hat\phi_0(\omega;
\gamma_j, z, {\rm sgn}(\gamma_{j+1}-\gamma_j))\\
\hat {W}_m^{\delta t}(\gamma, q; z, \omega) &=  \hat\psi_{\delta
t}(\omega, \gamma_q) \prod_{j=0}^{q-1} \hat\phi_{\delta t}(\omega;
\gamma_j, z, {\rm sgn}(\gamma_{j+1}-\gamma_j)).
\end{align}
We have that
\begin{align}
G^{\delta t}_{m}(x, y; z, \omega) - G_m(x, y; z, \omega) = & {1\over h}
\sum_{q=1}^\infty 2 ^{-q} \sum_{
\begin{matrix}
\gamma\in \Gamma_m : \gamma_0 = x, \gamma_{q} = y\\
\lvert\gamma_j - \gamma_{j-1}\lvert = 1 \forall j\ge1
\end{matrix}}
\left( \hat {W}_m^{\delta t} (\gamma, q; z, \omega) - \hat W_m (\gamma,
q; z, \omega)\right). \label{eq_udiff}
\end{align}
The integration over the contour in Fig.  \ref{fig_contour2} can
again be split into an integration over the countour ${\mathcal
C}_-$ and an integration over ${\mathcal C}_+$. The integral over
${\mathcal C}_+$ can be bounded from above thanks to Lemma
\ref{lemma_cplus}. Furthermore,  we have that
\begin{align}
&\bigg\lvert \int_{{\mathcal C}_-} \bigg( G^{\delta t}_{m}(x, y;
z, \omega) - G_m(x, y; z, \omega) \bigg)
e^{i\omega t}  {d\omega\over 2\pi} \bigg\lvert \notag \\
&\leq c h^{-1} \sqrt{q_{\rm max}} \max_{
\begin{matrix}
q, \gamma\in \Gamma_m : \gamma_0 = x, \gamma_{q} = y\\
\lvert\gamma_j - \gamma_{j-1}\lvert = 1 \forall j\ge1
\end{matrix}}
\bigg \lvert \int_{{\mathcal C}_-} \bigg( \hat {W}_m^{\delta t}
(\gamma, q; z, \omega) - \hat W_m (\gamma, q; z, \omega) \bigg) e^{i\omega
t}  {d\omega\over 2\pi}  \bigg\lvert. \notag
\\
& \leq c h^2 \label{eq_udiff2}
\end{align}

To bound the time derivative, we have to consider
\begin{align}
&\bigg\lvert \int_{{\mathcal C}_-} \bigg( {e^{i\omega\delta t} - 1
\over \delta t} G^{\delta t}_{m}(x, y; z, \omega) - i \omega G_m(x, y;
z, \omega) \bigg)
e^{i\omega t}  {d\omega\over 2\pi} \bigg\lvert \notag \\
\end{align}
But, since $\delta t = O(h^2)$, also this difference is $O(h^2)$.

\section{Convergence in the Ito Representation \label{SecItoRepresentation}}

The convergence Theorem \ref{kernel_convergence} admits a second formulation.
Let us introduce the function
\begin{equation}
\phi_m(x) = h_{xm} \sum_{x'\in X_m, x'\leq x} a(x').
\end{equation}
We have that
\begin{align}
e^{ip\phi_m(x)} \LLL(x, x') & e^{-ip\phi_m(x')} = \left(
{\sigma(x)^2\over 2 h_{xm}^2} + {\mu(x)\over 2 h_{xm}} \right) e^{-
i {h_{xm} a(x+h_{xm})} p } \delta_{x', x+h_{xm}} \notag \\
& + \left( {\sigma(x)^2\over 2 h_{xm}^2} - {\mu(x)\over 2 h_{xm}}
\right) e^{i h_{xm} a(x) p } \delta_{x', x-h_{xm}}
-{\sigma(x)^2\over h_{xm}^2} \delta_{x' x} \notag \\
& = {1\over2} \sigma(x, p, h_{xm})^2 \Delta^m(x, x') + \mu(x, p,
h_{xm}) \nabla^m(x, x')  \notag \\
&+\bigg[- i {\sin( h_{xm} p a(x)) \over h_{xm}}
\mu(x) + {1\over2} (\cos (h_{xm} a(x) p)-1)\bigg] \delta_{x x'}\notag\\
& + \left( {\sigma(x)^2\over 2} + {h_{xm} \mu(x)\over 2} \right)
 \bigg( h_{xm}^{-1} \nabla_x^+  e^{- i {h_{xm} a(x)} p} \bigg)
 \bigg( \delta_{x x'} + h_{xm} \nabla_x^+(x, x')  \bigg)\notag \\
\label{eq_fourierl}
\end{align}
where $\sigma(x, p, h_{xm})^2$ and $\mu(x, p, h_{xm})$ are defined
in (\ref{eq_sigma2}) and (\ref{eq_mu2}) and
\begin{equation}
\nabla^{+}_x f(x) = {f(x+{h_{xm}})-f(x)\over {h_{xm}}}.
\end{equation}
Hence, we have that
\begin{align}
\tilde\LLL^{m}(x, x'; p) & =  e^{ip\phi_m(x)} L^{m}(x, x'; p) e^{-ip\phi_m(x')} \notag \\
\label{eq_fourierl}
\end{align}
where
\begin{align}
L^{m}(x, x'; p) & =  \LLL(x, x') + \zeta_m(x, p, h_{xm}) \delta_{x x'}  + h_{xm} r_m(x, x'; p).
\label{eq_fourierl}
\end{align}
Here
\begin{align}
&\zeta_m(x, p, h_{xm}) =  i p b(x) +
 {\sigma(x)^2\over 2} \bigg( h_{xm}^{-1} \nabla_x^+  e^{- i {h_{xm} a(x)} p} \bigg),
\notag \\
&r_m(x, x'; p) = \bigg( h_{xm}^{-1} \nabla_x^+
e^{- i {h_{xm} a(x)} p} \bigg) \bigg[   \left( {\sigma(x)^2\over 2} + {h_{xm}
\mu(x)\over 2} \right) \nabla_x^+(x, x') + {\mu(x)\over 2} \delta_{x x'} \bigg] e^{ip(\phi_m(x') - ip\phi_m(x))}.
\notag \\
\label{eq_fourierl}
\end{align}
and
\begin{equation}
\nabla^{m+}_x f(x) = {f(x+{h_{xm}})-f(x)\over {h_{xm}}}.
\end{equation}

The operator $\tilde L^{m}(x, x'; p)$ is equivalent to
$\tilde \LLL^{m}(x, x'; p)$ up to a non-singular linear transformation.
We say that $\tilde L^{m}(x, x'; p)$ is the {\it Fourier transformed
generator in the Ito representation}. Also notice that the weak limit of this operator
as $m\to\infty$ is equal to
\begin{align}
\lim_{m\to\infty} L^{m}(x, x'; p)  =  {\sigma(x)^2\over 2} {\partial^2\over \partial x^2} +
\mu(x) {\partial\over \partial x} + i p b(x) - {i p \sigma(x)^2\over 2}{\partial a(x)\over \partial x}.
\end{align}

The Fourier transformed kernel is given by
\begin{eqnarray}
\hat u_m(x, x'; p) = e^{ip(\phi_m(x)-\phi_m(x'))} \hat U_m(x, x'; p)\notag \\
\end{eqnarray}
where
\begin{eqnarray}
\hat U_m(x, x'; p) = h_{xm}^{-1} e^{ t L^{m} }(x, x').\notag \\
\end{eqnarray}
The joint kernel is thus given by the following formula whose continuum analog
was found in \cite{Girsanov},
\cite{CameronMartin}, \cite{Feynman1948} and  \cite{Ito1951}:
\begin{eqnarray}
u(x, I; x', I'; t) = \int_{-\infty}^\infty {dp\over 2\pi} e^{ip [ I' - I -\phi_m(x') + \phi_m(x)]}
\hat U_m(x, x'; p).\notag \\
\end{eqnarray}

\begin{theorem}{\label{kernel_convergence2}} Under the same assumptions of
Theorem \ref{kernel_convergence}, for all $K>0$ there is a constant $c(K)$ such that
\begin{equation} \lvert\lvert \hat U_{m}^\phi(t) - \hat U_{m'}^\phi(t)
\lvert\lvert_{m, K, \hat {\mathcal L}} \leq c(K) h_{xm}^{\gamma}.
\end{equation}
for all $m'>m\ge m_0$. A similar bound also holds for the kernels
obtained with a fully explicit Euler scheme, i.e.
\begin{equation}
\hat U^{\delta t}_{m}(x, x'; p, t) = {1\over h_{x}} \left( 1 +
{\delta t}_{m} L^{m}(p) \right)^{\left[t\over\delta
t_{m}\right]}(x, x'),
\end{equation}
where $\delta t_{m} $ is so small that
\begin{equation}
\min_{x, p\in X_m\times \hat Y_n} 1 + \delta t_{m} \hat {\mathcal L}^{m}(x, x; p)
> 0.
\end{equation}
In this case, there is a constant $c(K)$ such that
\begin{equation} \lvert\lvert U_{m}(t) - U^{\delta t}_{m}(t)
\lvert\lvert_{m, K, \hat {\mathcal L}} \leq c(K) h_{xm}^2
\end{equation}
for all $m'>m\ge m_0$.
\end{theorem}

It is possible to retrace the argument in the previous section except for replacing
\begin{equation}
\kappa(x, p, h_{xm}) \rightarrow \zeta_m(x, p, h_{xm}) + h_{xm} r_m(x, p, h_{xm})
\end{equation}
where $\kappa(x, p, h_{xm})$ is the function in (\ref{eq_kappa2}). All arguments go through
unaltered under the same conditions. As a consequence we conclude that the kernels
\begin{eqnarray}
{1\over h_{xm}} \exp\bigg(t \big(\LLL^m + \zeta_m(z) + h_{xm} r_m(z) \big)\bigg)(x, x')
\end{eqnarray}
converge in graph norm, uniformly on discs $z\in {\Bbb C} : \lvert z \lvert < K$, for any $K>0$

\section{Other Abelian Processes\label{SupProcess}}

In this Section, we give to examples of Abelian processes which
emerge from applications and are not stochastic integrals.

\subsection{The Sup Process}

Consider the sup process
\begin{equation}
y_t = \sup_{\in [0,t]} x_s.
\end{equation}
the sup is always attained as an element of the
underlying space $X_m$, it is natural in this case
to restrict the attention to the case where $n=m$ and
$Y_m = X_m$. The joint generator is given by
\begin{align}
\tilde {\mathcal L}^{m}  (x, y; x', y') =& {\mathcal L}^{m}  (x; x')
{\mathcal A}^{m}  (x, y; x', y') \notag
\\
\end{align}
where $x, y, x', y' \in X_m$ and we set
\begin{align}
{\mathcal A}^{m}  (x, y; x', y') =
\begin{cases}
\delta_{y y'} \;\;\;\;\;\;{\rm if}\;\;\;\;x' < y \\
\delta_{x' y'} \;\;\;\;\;{\rm if}\;\;\;\;x' \ge y.
\end{cases}
\end{align}
Consider the matrix
\begin{align}
{\mathcal V}_{m}  (x, y; x', y') = \delta(x - x') 1(y' \ge y)
\end{align}
and its inverse
\begin{align}
{\mathcal V}_{m}^{-1}  (x, y; x', y') = \delta(x - x')
\big(\delta(y'-y) - \delta(y'-y - h_{xm})\big).
\end{align}
Consider the one parameter family of operators $\hat {\mathcal
L}^{m} (y) $ such that
\begin{align}
\big( {\mathcal V}_{m}^{-1} \tilde {\mathcal L}^{m} {\mathcal V}_{m}
\big) (x, y; x', y') =  \hat {\mathcal L}^{m} (x, x'; y) \delta_{y
y'}.
\end{align}
We have that
\begin{align}
\hat {\mathcal L}^{m} (x, x'; y) = {\mathcal L}^{m}  (x; x') 1(x'
\leq y).
\end{align}
Hence, $\hat {\mathcal L}^{m} (x, x'; y)$ is the Markov generator
of the underlying process with absorption in the interval $[y, L_x]$.

Let us notice that the kernel can be obtained as follows:
\begin{equation}
u_{m}(t) = {\mathcal V}_{m} \tilde u_{m}(t) {\mathcal V}_{m}^{-1}
\;\;\;\;{\rm where}\;\;\;\; \tilde u_{m}(x, y; x', y'; t) =
\delta(y- y') \exp\left( t \tilde {\mathcal L}^{m}(y) \right)(x,
x').
\end{equation}
A more explicit way of expressing the joint kernel is
\begin{equation}
u_{m}(x, y; x', y'; t) =
\delta(y - y') \hat u_{m}((x, x'; y, t)
+
1(y' > y)\bigg(
\hat u_{m}(x, x'; y', t)
-
\hat u_{m}(x, x'; y'-h_{xm}, t)\bigg)
\end{equation}

Convergence in this case can be established along the same lines as done for
stochastic integrals. The situation is simpler in that only the consideration
of the kernel itself, i.e. the $p=z=0$ case with the notations in Sections 2 and 3,
is needed. The additional complication is that we need to consider absorbing
boundary conditions. This implies a few marginal changes to the derivation above,
as when a path arrives to an absorbing lattice point, it stays constant
thereafter. Since this applies both to paths and decorating paths in a finer lattice,
upon arriving to an absorption point the dynamics is trivial in either case
and the final bounds given still hold.

\subsection{Discrete Time Processes\label{DiscreteTime}}

This section is based on work in collaboration with Manlio Trovato
\cite{ATrovato2} and Paul Jones, see \cite{AJones}.

An important class of path-dependent options requires computing the
joint distribution of the underlying lattice process and of a
discrete sum of the following form:
\begin{equation}
y_t = \sum_{i=1}^N \psi(x_{t_{i-1}}, x_{t_{i}};
t_i) \label{eq_defi2}
\end{equation}
where $N$ is an integer, $t_i = i \Delta T$ and $T = N \Delta T$.
Consider the elementary propagator
\begin{equation}
U_{m}(x_1, x_2) = {1\over h_{xm}} e^{(\Delta T) \LLL} (x_1, x_2).
\end{equation}
To find the joint transition probability, one can again discretize
the variable $y_t$ in the lattice $Y_n = h_{yn} \ZZZ \cap [-L_y, L_y]$.
As opposed to lifting the generator as done above for the other cases,
here we lift the elementary propagator itself and form the joint propagator
\begin{equation}
\tilde U_{mn}(x_1, y_1; x_2, y_2) = U_{m}(x_1, x_2) \delta\big(y_1 -
y_2 + [ \psi(x_1, x_2) h_{yn}^{-1}]\big).
 \label{eq_defu}
\end{equation}
This lifted operator can be block-diagonalized by means of a partial
Fourier transform. Consider the Fourier transform operator $\hat U_{m}(p)$ of
matrix elements
\begin{equation}
\hat U_{m}(x_1 , x_2; p) = \lim_{n\to\infty}
U_{mn}(x_1, y_1; x_2, y_2) e^{ -i p (y_2 - y_1)} =
U_{m}(x_1, x_2) e^{ -i p \psi(x_1,x_2)}.
 \label{eq_defup}
\end{equation}
Then we have that
\begin{equation}
\lim_{n\to\infty} \big(\tilde U_{mn}^N\big) (x_1, y_1; x_2, y_2) =
\int_{-\infty}^\infty {dp\over 2\pi} e^{ip(y_2-y_1)} \big(
\hat U_m^N \big)(x_1, x_2; p)
\label{discretefk1}
\end{equation}
Convergence in the graph-norm in this case descends directly from the convergence
of the one-period kernel $U_{m}(x_1, x_2)$.

\section{Conclusions}

We obtained bounds on convergence rates for explicit discretization
schemes to the Fourier transform of joint
kernels of one-dimensional diffusion equations with
continuous coefficients and a stochastic integral.
We consider both semi-discrete triangulations with continuous time and
explicit Euler schemes with
time step small enough for the method to be stable. The proof is
constructive and based on a new technique of path conditioning for
Markov chains and a renormalization group argument. Convergence
rates depend on the degree of smoothness and H\"older
differentiability of the coefficients. The method also applies
to a more general class of path dependent processes we call Abelian. Examples
of Abelian processes beside stochastic integrals are the sup process and
discrete time summations.

\bibliographystyle{giwi}
\bibliography{joint}

\end{document}